\documentclass[11pt, reqno]{amsart}
\usepackage{amsmath, amsthm, amscd, amsfonts, amssymb, graphicx, color}
\usepackage[bookmarksnumbered,plainpages]{hyperref}

\newtheorem{theorem}{Theorem}[section]
\newtheorem{lemma}[theorem]{Lemma}
\newtheorem{proposition}[theorem]{Proposition}
\newtheorem{corollary}[theorem]{Corollary}
\theoremstyle{definition}

\theoremstyle{remark}
\newtheorem{remark}[theorem]{Remark}
\numberwithin{equation}{section}

\begin{document}
\title{Constructions preserving  $n$-weak amenability of Banach algebras}
\author{A. Jabbari}
\address{Department of Mathematics, Ferdowsi University of Mashhad\\
P. O. Box 1159, Mashhad 91775, Iran.}\email{shahzadeh@math.um.ac.ir}
\author{M.\ S.\ Moslehian}
\address{Department of Mathematics, Ferdowsi University of Mashhad\\
P. O. Box 1159, Mashhad 91775, Iran; \newline Centre of Excellence
in Analysis on Algebraic Structures (CEAAS), Ferdowsi University of
Mashhad, Iran.} \email{moslehian@ferdowsi.um.ac.ir and
moslehian@ams.org}

\author{H.\ R.\ E.\ Vishki}
\address{Department of Mathematics, Ferdowsi University of Mashhad\\
P. O. Box 1159, Mashhad 91775, Iran; \newline Centre of Excellence
in Analysis on Algebraic Structures (CEAAS), Ferdowsi University of
Mashhad, Iran.} \email{vishki@ferdowsi.um.ac.ir}


\begin{abstract}
A surjective bounded homomorphism fails to preserve $n$-weak
amenability, in general. We however show that it preserves the
property if the involved homomorphism enjoys a right inverse. We
examine this fact for certain homomorphisms on several Banach
algebras.
\end{abstract}
\maketitle


\section{Introduction}

Let ${\mathcal A}$ be a Banach algebra and let ${\mathcal X}$ be a
Banach ${\mathcal A}$-bimodule or simply an ${\mathcal A}$-bimodule. A
bounded linear mapping $D: {\mathcal A} \longrightarrow {\mathcal
X}$ is said to be a derivation if $D(ab)=a\cdot D(b)+D(a)\cdot b$
for all $a, b\in {\mathcal A}$. A derivation $D : {\mathcal
A}\longrightarrow {\mathcal X}$ is said to be inner if there exists
$x\in {\mathcal X}$ such that $D(a)=a\cdot x-x\cdot a$ for all $a\in
{\mathcal A}$; in this case we say the inner derivation $D$ is
implemented by $x$.\\
For a Banach ${\mathcal A}$-module ${\mathcal X}$, the dual
${\mathcal X}^*$ of ${\mathcal X}$ equipped with the actions
$(f.a)(x)=f(ax)$ and $(a.f)(x)=f(xa)$ is a Banach ${\mathcal
A}$-module. Similarly, the $n$-th dual ${\mathcal X}^{(n)}$ of
${\mathcal X}$ is a Banach ${\mathcal A}$-module. In particular,
${\mathcal A}^{(n)}$ is a Banach ${\mathcal A}$-module.

 For terminology and background materials we follow \cite{D, HR, RUN}.

The notion of weak amenability was first introduced by Bade, Curtis
and Dales \cite{BCD} for commutative Banach algebras and then by
Johnson in \cite{J} for arbitrary Banach algebras. A Banach algebra
$\mathcal A$ is called weakly amenable if every bounded derivation
from $\mathcal A$ into the dual Banach module ${\mathcal A}^*$ is
inner. The concept of $n$-weak amenability was initiated by Dales,
Ghahramani and Gronb\ae k \cite{DGG}. A Banach algebra ${\mathcal
A}$ is said to be $n$-weakly amenable ($n\in\mathbb{N}$) if every
bounded derivation from ${\mathcal A}$ into ${\mathcal A}^{(n)}$ is
inner. Trivially, $1$-weak amenability is nothing else than weak
amenability. In \cite{DGG} the authors   presented many substantial
properties of this variety of Banach algebras. For instance, they
showed that $C^*$-algebras are $n$-weakly amenable for all
$n\in{\mathbb N}$ as well as that $L^1(G)$ is $(2n-1)$-weakly amenable
for all $n\in {\mathbb N}$ and for any locally compact group $G$.

One of the basic tools for making a new amenable Banach algebra from
the old one is the fact that: a continuous homomorphic image of an
amenable Banach algebra is again amenable. The same is also easily
valid for weak amenability in the commutative case, but is false in
general. An interesting counterexample is given by Gronb\ae k
\cite{Gro1}. (Indeed, the tensor algebra $X\widehat\otimes X^*$ is
weakly amenable for any Banach space $X$. However, the quotient
${\mathcal N}(X)$, the nuclear operators,
 is weakly amenable if and only if the kernel has
dimension  less than or equal to $1$. Taking $X=E\oplus E$  with a
Banach space $E$ failing the approximation property, we get an
example of a weak amenable Banach algebra whose quotient by a closed
ideal is not weakly amenable.) The interested reader is referred to
\cite{Gro2} for a sufficient condition for weak amenability of the
homomorphic image.

In this paper we give a sufficient condition for the $n$-weak
amenability of the homomorphic image. More precisely, we show that the
homomorphic image of an $n$-weakly amenable Banach algebra is
$n$-weakly amenable if  the involved homomorphism enjoys a right
inverse. We use this fact for obtaining some $n$- weakly amenable
Banach algebras from old ones. Moreover, our method in turn provides
a unified approach from which  some results of Ghahramani and Laali
\cite{GL} and also Lau and Loy \cite{LL} can be directly derived.


\section{Right invertible homomorphisms  preserve $n$-weak amenability}

For a Banach algebra $\mathcal A$, we denote the natural ${\mathcal
A}$-module actions on ${\mathcal A}^*$ by
\[\langle a\cdot a^*, b\rangle=\langle a^*, ba\rangle, \  \ \ \ \ \langle a^*\cdot a, b\rangle=\langle a^*, ab\rangle
\ \ \ \ \ (a, b\in {\mathcal A}, a^* \in {\mathcal A}^*).\] In the
case when $\Theta : {\mathcal A}\longrightarrow {\mathcal B}$ is a
bounded homomorphism, ${\mathcal B}^{(n)}$ can be regarded as an
${\mathcal A}$-module under the module actions

\[a\cdot b^{(n)}=\Theta(a)\cdot b^{(n)}, \ \ \
b^{(n)}\cdot a= b^{(n)}\cdot \Theta(a)\ \ \  \ (a\in {\mathcal A},
b^{(n)}\in {\mathcal B}^{(n)}).\] A direct verification reveals that
the adjoint mappings $\Theta^{*} : {\mathcal B}^{*}\longrightarrow
{\mathcal A}^{*}$ and $\Theta^{**} : {\mathcal
A}^{**}\longrightarrow {\mathcal B}^{**}$ are ${\mathcal B}$-module
morphisms. The same fact is true for the higher adjoint mappings
$\Theta^{(2n-1)} : {\mathcal B}^{(2n-1)}\longrightarrow {\mathcal
A}^{(2n-1)}$ and $\Theta^{(2n)} : {\mathcal A}^{(2n)}\longrightarrow
{\mathcal A}^{(2n)}$.\\

We commence with the next lemma which plays a key role in the
sequel.
\begin{lemma} \label{lemma} Let ${\mathcal A}$ and ${\mathcal B}$ be  Banach
algebras and let $\Theta : {\mathcal A}\longrightarrow {\mathcal B}$
and $\Phi : {\mathcal B}\longrightarrow {\mathcal A}$ be bounded
homomorphisms such that $\Theta\circ\Phi=I_{\mathcal B}$ and
$n\in\mathbb{N}$.

$(i)$ If $D : {\mathcal B}\longrightarrow {\mathcal B}^{(2n-1)}$ is a
derivation then so is $D_1=(\Theta^{(2n-1)}\circ D\circ \Theta) :
{\mathcal A}\longrightarrow {\mathcal A}^{(2n-1)}$. Moreover, if
$D_1$ is inner then so is $D$.

$(ii)$ If $D : {\mathcal B}\longrightarrow {\mathcal B}^{(2n)}$ is a
derivation then so is $D_2=(\Phi^{(2n)}\circ D\circ \Theta) :
{\mathcal A}\longrightarrow {\mathcal A}^{(2n)}$. Moreover, if $D_2$
is inner then so is $D$.

$(iii)$ If ${\mathcal A}$ is $n$-weakly amenable then so is
${\mathcal B}$.
\end{lemma}
\begin{proof}
$(i)$ As we have mentioned before,  $\Theta^{(2n-1)}: {\mathcal
B}^{(2n-1)}\longrightarrow {\mathcal A}^{(2n-1)}$ is an ${\mathcal
A}$-module morphism. Now let $D : {\mathcal B}\longrightarrow
{\mathcal B}^{(2n-1)}$ be a derivation. Then for every $a, c\in
{\mathcal A},$
\begin{eqnarray*} D_1(ac)=(\Theta^{(2n-1)}\circ D\circ\Theta)(ac)
&=&\Theta^{(2n-1)}\circ
D(\Theta(a)\Theta(c))\\&=&\Theta^{(2n-1)}(a\cdot
D(\Theta(c))+D(\Theta(a))\cdot
c)\\&=&a\cdot\Theta^{(2n-1)}(D(\Theta(c)))+\Theta^{(2n-1)}(D(\Theta(a)))\cdot
c\\&=&a\cdot D_1(c)+D_1(a)\cdot c,\end{eqnarray*} which means that
$D_1 : {\mathcal A}\longrightarrow {\mathcal A}^{(2n-1)}$ is a
derivation.

Now assume that $D_1$ is inner which is implemented by $F\in
{\mathcal A}^{(2n-1)}$. As $\Theta\circ\Phi=I_{\mathcal B}$, we have
$\Theta^{(2n-2)}\circ\Phi^{(2n-2)}=I_{{\mathcal B}^{2n-2}}$, and so
for every $b\in {\mathcal B}$ and $G\in B^{2n-2}$ (with ${\mathcal
B}^{(0)}={\mathcal B}$ in mind) we have
\begin{eqnarray*} \langle D(b), G\rangle&=&  \langle D(\Theta\circ\Phi(b)),
(\Theta^{(2n-2)}\circ\Phi^{(2n-2)})(G)\rangle\\&=&\langle D_1(\Phi(b)), \Phi^{(2n-2)}(G)\rangle\\
&=&\langle \Phi(b)\cdot F-F\cdot\Phi(b), \Phi^{(2n-2)}(G)\rangle\\&=&\langle \Phi^{(2n-1)}(b\cdot F-F\cdot b), G\rangle\\
&=&\langle b\cdot\Phi^{(2n-1)}(F)-\Phi^{(2n-1)}(F)\cdot b, G\rangle\,.\\
\end{eqnarray*}
(Note that in the penultimate identity we used the fact that
$\Phi^{(2n-1)} : {\mathcal A}^{(2n-1)}\longrightarrow {\mathcal
B}^{(2n-1)}$ is actually a ${\mathcal B}$-module morphism.)
Therefore, $D$ is an inner derivation implemented by
$\Phi^{(2n-1)}(F)$. The proof of $(ii)$ is similar to the  proof of
$(i)$. Part $(iii)$ follows trivially from $(i)$ and $(ii)$.
\end{proof}


\begin{remark} It is worthwhile to mention that  the hypothesis imposed in
Lemma~\ref{lemma} is equivalent to the  fact that $\mathcal A$ is
the topological direct sum of Banach spaces $\ker(\Theta)$ and $B$.
Indeed,  direct verification reveals that the mapping $\omega:
{\mathcal A} \to \ker(\Theta) \times {\mathcal B}$ defined by
$\omega(a)=\left(a-\Phi(\Theta(a)), \Theta(a)\right)$ is a
bicontinuous isomorphism whose inverse is $\varpi: \ker(\Theta)
\times {\mathcal B} \to {\mathcal A}$ defined by
$\varpi(a,b)=a+\Phi(b)$.
\end{remark}

As an  application of  Lemma~\ref{lemma} we have the following result
which shows that $n$-weak amenability inherits by those closed
subalgebras which are  complemented by a closed ideal.


\begin{theorem}\label{M}
 Let  ${\mathcal A}$ be a Banach algebra such that  ${\mathcal A}={\mathcal B}\oplus
{\mathcal I}$ for some closed ideal ${\mathcal I}$ and closed
subalgebra ${\mathcal B}$. Then $n$-weak amenability of ${\mathcal
A}$ implies that of ${\mathcal B}$.
\end{theorem}

\begin{proof}
Let $\Theta : {\mathcal A}\longrightarrow {\mathcal B}$ be the
natural projection onto $B$ and let $\Phi : {\mathcal
B}\longrightarrow {\mathcal A}$ be the natural injection into
${\mathcal A}$. Trivially $\Theta$ and $\Phi$ are (bounded)
homomorphisms with $\Theta\circ\Phi=I_{\mathcal B}$. Now the
conclusion follows from Lemma~\ref{lemma}.
\end{proof}

The next result is a restatement of the above result in the
language of quotient algebras and split short exact sequences.
Recall that a short exact sequence $0\rightarrow {\mathcal
C}\rightarrow {\mathcal A}\rightarrow {\mathcal B}\rightarrow 0$ in
the category of Banach algebras and bounded homomorphisms splits if
the morphism ${\mathcal A}\rightarrow {\mathcal B}$ enjoys a right
inverse morphism; cf. \cite{H}.


\begin{corollary} If  ${\mathcal A}$ is $n$-weakly amenable and ${\mathcal I}$ is a closed ideal of ${\mathcal A}$
then the quotient algebra ${\mathcal A}/{\mathcal I}$ is $n$-weakly
amenable provided the natural short exact sequence $0\rightarrow
{\mathcal A}\rightarrow {\mathcal A}\rightarrow {\mathcal
A}/{\mathcal I}\rightarrow 0$ splits.
\end{corollary}

Let ${\mathcal X}$ be a  Banach ${\mathcal A}$-module. Then the
$l^1-$direct sum ${\mathcal A}\oplus {\mathcal X}$ is a Banach
algebra under
\[(a, x)\cdot (b,y)=(ab, ay+xb) \ \ \ \ \ (a, b\in {\mathcal A}, x, y\in {\mathcal X})\,,\]
which is known as a module extension Banach algebra. Some properties
of algebras of this form have been studied in~\cite{DGG}. In
particular, the $n$-weakly amenability of this kind of Banach
algebras has been extensively investigated in \cite[Theorems 2.1 and
2.2]{Z}. Since ${\mathcal X}$ and ${\mathcal A}$ are ideal and
closed subalgebra of ${\mathcal A}\oplus {\mathcal X}$,
respectively, as a consequence of Theorem \ref{M} we have


\begin{corollary}
Let ${\mathcal X}$ be a Banach ${\mathcal A}$-module. If the
module extension Banach algebra ${\mathcal A}\oplus {\mathcal X}$ is
$n$-weakly amenable then so is ${\mathcal A}$.
\end{corollary}


\section{$n$-weak amenability of second dual}

Recall that for a Banach algebra ${\mathcal A}$, the  second dual
${\mathcal A}^{**}$ of ${\mathcal A}$ can be made into a Banach
algebra supplied with either the first Arens product $\Box$ or the second
Arens product $\lozenge$. For $a\in {\mathcal A}$, $a^*\in {\mathcal
A}^*$ and $a^{**}, b^{**}\in {\mathcal A}^{**}$, the elements
$a^{**}\cdot a^*$ and $a^*\cdot a^{**}$ of ${\mathcal A}^{*}$ are
defined by the formulas
\[\langle a^{**}\cdot a^*, a\rangle=\langle a^{**}, a^*\cdot a\rangle, \ \ \ \ \
\langle  a^{*}\cdot a^{**}, a\rangle=\langle a^{**}, a\cdot
a^*\rangle.\] Next, $a^{**}\Box \ b^{**}$ and $a^{**}\lozenge \
b^{**}$ are defined in ${\mathcal A}^{**}$ by the formulas
\[\langle a^{**}\Box \ b^{**}, a^*\rangle=\langle a^{**}, b^{**}\cdot a^*\rangle, \ \
\ \ \ \langle  a^{**}\lozenge  \ b^{**}, a^*\rangle=\langle b^{**},
a^*\cdot a^{**}\rangle.\] Ample information about Arens products may
be found in  Arens' original paper \cite{AR} (see also \cite{D}).

For a Banach space ${\mathcal X}$, when there is no risk of
confusion, we usually identify ${\mathcal X}$ $(x\in {\mathcal X})$
with its canonical image in ${\mathcal X}^{**}$.

A Banach algebra ${\mathcal A}$ is said to be a dual Banach algebra
if there is a closed submodule ${\mathcal A}_*$ of ${\mathcal A}^*$
for which ${\mathcal A}=({\mathcal A}_*)^*$ (${\mathcal A}_*$ is
called the predual of ${\mathcal A}$). As a consequence of
Theorem~\ref{M} we shall see that in the case when ${\mathcal A}$
is a dual Banach algebra the $n$-weak amenability of ${\mathcal
A}^{**}$ implies that of ${\mathcal A}$; see \cite[Theorem 2.2]{GL}
for $n=1$.


\begin{proposition} \label{M1}
For a dual Banach algebra ${\mathcal A},$ if ${\mathcal A}^{**}$
(supplied with either of the Arens products) is $n$-weakly amenable,
then so is ${\mathcal A}$.
\end{proposition}
\begin{proof}
We present a proof for the first Arens product $\Box$. Let ${\mathcal
A}$ be a dual Banach algebra with respect to ${\mathcal A}_*$.
Suppose that $\Theta : {\mathcal A}^{**}\longrightarrow {\mathcal
A}$ is the adjoint of the canonical embedding $J : {\mathcal
A}_*\longrightarrow ({\mathcal A}_*)^{**}={\mathcal A}^*$ and let
$\Phi : {\mathcal A}\longrightarrow {\mathcal A}^{**}$  be also the
canonical embedding. Then trivially $\Theta\circ\Phi=I_{\mathcal A}$
and  also $\Phi : {\mathcal A}\longrightarrow ({\mathcal A}^{**},
\Box)$ is a  homomorphism. Moreover, $\Theta : ({\mathcal A}^{**},
\Box)\longrightarrow {\mathcal A}$ is a homomorphism. Indeed, let
$a_*\in {\mathcal A}_*$, $a^{**}, b^{**}\in {\mathcal A}^{**}$ and
let $\{a_\alpha\}$, $\{b_\beta\}$ be two nets in ${\mathcal A}$ that
converge  to $a^{**}$ and $b^{**}$ in the $w^*-$topology on
${\mathcal A}^{**}$, respectively. Then
\begin{eqnarray*}\langle \Theta(a^{**}\Box \ b^{**}), a_*\rangle&=&\langle a^{**}\Box \ b^{**}, J(a_*)\rangle
=\lim_\alpha\lim_\beta\langle a_\alpha b_\beta,
J(a_*)\rangle\\&=&\lim_\alpha\lim_\beta\langle
\Theta(a_\alpha)\Theta(b_\beta), a_*\rangle=\lim_\alpha\langle
J(a_*\cdot\Theta(a_\alpha)), b^{**}\rangle\\&=&\lim_\alpha\langle
a_\alpha, J(\Theta(b^{**})\cdot a_*)\rangle=\langle a^{**},
J(\Theta(b^{**})\cdot a_*)
\rangle\\&=&\langle\Theta(a^{**})\Theta(b^{**}),
a_*\rangle.\end{eqnarray*} Using Theorem~\ref{M} we conclude the
result.
\end{proof}


For a Banach algebra ${\mathcal A}$ let ${\mathcal A}^{op}$ be the
Banach algebra endowed with the  reversed product of ${\mathcal A}$,
whose  underlying Banach space is ${\mathcal A}$ itself. The
$n$-weak amenability of ${\mathcal A}$ is trivially equivalent to
that of ${\mathcal A}^{op}$. As another application of
Theorem~\ref{M} we have the following result; see \cite[Theorem
2.3]{GL} for $n=1$.


\begin{proposition}
Let ${\mathcal A}$ be a Banach algebra admitting a continuous
antihomomorphism $\lambda$ such that $\lambda^2=I_{\mathcal A}$.
Then $({\mathcal A}^{**}, \Box)$ is $n$-weakly amenable if and only
if $({\mathcal A}^{**}, \lozenge)$ is $n$-weakly amenable. A similar
conclusion holds if $\lambda$ is a continuous involution on
${\mathcal A}$.
\end{proposition}
\begin{proof} As $\lambda^2=I_{\mathcal A}$, for the second adjoint
$\lambda^{**}$ of $\lambda$ we have $(\lambda^{**})^2=I_{{\mathcal
A}^{**}}.$  Using a limit process similar to what we used in the
proof of Proposition~\ref{M1}, one can show that $\lambda^{**} :
({\mathcal A}^{**}, \Box)\longrightarrow ({\mathcal A}^{**},
\lozenge)^{op}$ and $\lambda^{**} :  ({\mathcal A}^{**},
\lozenge)^{op}\longrightarrow ({\mathcal A}^{**}, \Box)$ are
homomorphisms. The conclusion follows from Theorem~\ref{M}. In the
case when $\lambda$ is an involution on ${\mathcal A}$  a similar
proof may apply.
\end{proof}


\section{Introverted subspaces of ${\mathcal A}^*$}

Let ${\mathcal A}$ be a Banach algebra. A closed subspace ${\mathcal
X}$ of ${\mathcal A}^*$ is said to be left invariant if ${\mathcal
X}\cdot{\mathcal A}\subseteq {\mathcal X}$ (or equivalently,
$\mathcal X$ is a right Banach ${\mathcal A}-$submodule of
${\mathcal A}^*$). A left invariant subspace ${\mathcal X}$ of
${\mathcal A}^*$ is said to be left introverted if ${\mathcal
X}^*\cdot {\mathcal X}\subseteq {\mathcal X}$. The dual ${\mathcal
X}^*$ of a left introverted subspace ${\mathcal X}$ of ${\mathcal
A}^*$ has a natural  (first Arens type) multiplication $\Box$ on it
in the same manner as ${\mathcal A}^{**}$ does. Right invariance and
right introversion for a subspace ${\mathcal X}$ of ${\mathcal A}^*$
can be defined similarly. In the case when ${\mathcal X}$ is right
introverted then ${\mathcal X}^*$ can be supplied with a natural
(second Arens type) multiplication $\lozenge$. A trivial example of
a (left and right) introverted subspace is ${\mathcal A}^*$ itself.
More illuminating examples of introverted subspaces of ${\mathcal
A}^*$ are $WAP({\mathcal A})$ (=weakly almost periodic elements of
${\mathcal A}^*$) and $AP({\mathcal A})$ (=almost periodic elements
of ${\mathcal A}^*$). It can be shown that every $w^*-$closed
invariant subspace of ${\mathcal A}^*$ is introverted. The same is
valid for every (norm) closed subspace of $WAP({\mathcal A})$; see
\cite[Lemma 1.2]{LL}. In the case when ${\mathcal A}$ enjoys a
bounded right  (left) approximate identity, the
Cohen-Hewitt Factorization Theorem \cite[Theorem 32.22]{HR} shows
that ${\mathcal A}^*\cdot{\mathcal A}$ (${\mathcal
A}\cdot {\mathcal A}^*$) is a closed subspace of ${\mathcal A}^*$.
Moreover, ${\mathcal A}^*\cdot{\mathcal A}$ is left
introverted while ${\mathcal A}\cdot{\mathcal A}^*$ is right introverted.\\

For every two left introverted subspaces  ${\mathcal Y}\subseteq
{\mathcal X}$ of ${\mathcal A}^*$ it can be readily verified that
the restriction map $P : {\mathcal X}^*\longrightarrow {\mathcal
Y}^*$ is a continuous homomorphism  onto ${\mathcal Y}^*$ whose
kernel  is the $w^*-$closed ideal ${\mathcal Y}^\bot=\{x^*\in
{\mathcal X}^* : \langle x^*, y\rangle=0  ~~{\rm for \ all}~~~ y\in
{\mathcal Y}\}$ of ${\mathcal X}^*$; see \cite[Lemma 1.1]{LL}. This
provides the direct sum decomposition ${\mathcal X}^*={\mathcal
Y}^*\oplus {\mathcal Y}^\bot$. A problem which is of interest in
this section is:  for two left introverted subspaces ${\mathcal
Y}\subseteq {\mathcal X}$ of ${\mathcal A}^*$, under what conditions
the $n$-weak amenability of ${\mathcal X}^*$ implies that of
${\mathcal Y}^*$? According to Theorem~\ref{M} it is the case if
${\mathcal Y}^*$ is isometrically isomorphic to a closed subalgebra
of ${\mathcal X}^*$. In this case we say the pair $({\mathcal Y}^*,
{\mathcal X}^*)$ is admissible. We summarize these observations in
the next result.
\begin{theorem}\label{M2}
Let ${\mathcal Y}, {\mathcal X}$ be two left introverted subspaces
of ${\mathcal A}^*$ such that ${\mathcal Y}\subseteq {\mathcal X}$.
If the pair $({\mathcal Y}^*, {\mathcal X}^*)$ is admissible then
$n$-weak amenability of ${\mathcal X}^*$ implies that of ${\mathcal
Y}^*$.
\end{theorem}
\begin{corollary}\label{M3}
Let $A$ have a right approximate identity bounded by $1$. Then the
$n$-weak amenability of ${\mathcal A}^{**}$ implies that of
$({\mathcal A}^*\cdot{\mathcal A})^*$.
\end{corollary}
\begin{proof} Using the latter result it is enough to show that the
pair $(({\mathcal A}^*\cdot{\mathcal A})^*, {\mathcal A}^{**})$ is
admissible. To see this we shall show that $({\mathcal
A}^*\cdot{\mathcal A})^*$ is isometrically isomorphic to the closed
subalgebra $e^{**}{\mathcal A}^{**}$ of ${\mathcal A}^{**}$, in
which $e^{**}$ is a $w^*-$cluster point of the involved bounded
approximate identity in ${\mathcal A}^{**}$ which in turn is a right
identity for ${\mathcal A}^{**}$. It can be readily verified that
the restriction map $\pi : e^{**}{\mathcal A}^{**}\longrightarrow
({\mathcal A}^*\cdot{\mathcal A})^*$ is a norm decreasing  (Banach
algebra) isomorphism, but not necessarily isometric. As the involved
right approximate identity  is bounded by $1$, $\pi$ actually is
an isometry. In fact, for every $a^{**}\in e^{**}{\mathcal A}^{**}$
we have trivially $e^{**}\Box \ a^{**}=a^{**}$ and so for each
$a^*\in {\mathcal A}^*$,
\begin{eqnarray*}
|\langle a^{**}, a^*\rangle|&=&|\langle e^{**}\Box \ a^{**},
a^*\rangle|\\
&=&\lim_\alpha|\langle \pi(a^{**}), a^*\cdot
e_{\alpha}\rangle|\\
&\leq& \lim_\alpha\|\pi(a^{**}) \| \|a^*\| \|e_{\alpha}\|\\
&\leq& \|\pi(a^{**}) \| \|a^*\|\,, \end{eqnarray*} which means that
$({\mathcal A}^*\cdot{\mathcal A})^*$ is isometrically isomorphic to
the closed subalgebra $e^{**}{\mathcal A}^{**}$ of ${\mathcal
A}^{**}$.
\end{proof}


Since  ${\mathcal A}=L^1(G)$, where $G$ is a locally compact group,
always enjoys an approximate identity of norm $1$, using the fact
that ${\mathcal A}^*\cdot{\mathcal A}=L^\infty(G)\cdot
L^1(G)=LUC(G)$, as a rapid application of the latter result we have
the following corollary. A  special version of it (for the case
$n=1$) is given in \cite[Proposition 4.14]{LL}.

\begin{corollary}
If $L^1(G)^{**}$ is $n$-weakly amenable then so is $LUC(G)^*$.
\end{corollary}


A result of Leptin {\cite{L}} asserts that if $G$ is amenable then
the Fourier algebra ${\mathcal A}=A(G)$ has a bounded (by $1$)
approximate identity (and vice versa). Also  in this case, ${\mathcal
A}^*\cdot{\mathcal A}=UC(\hat{G})^*$. Applying Corollary~\ref{M3} we
have the following result, whose special case $n=1$ is given in
\cite[Proposition 6.3]{LL}.

\begin{corollary}
If $G$ is an amenable group then $n$-weak amenability of
$A(G)^{**}$ implies the same for $UC(\hat{G})^*$.
\end{corollary}


\begin{remark} (i) Note that in the case when  $\mathcal A$
is a dual Banach algebra (with predual ${\mathcal A}_*$) then
trivially ${\mathcal A}_*$ is an introverted subspace of ${\mathcal
A}^*$ (which in turn supports the decomposition ${{\mathcal
A}}^{**}={\mathcal A}\oplus {A_*}^\bot$). As $({\mathcal
A}_*)^*={\mathcal A}$ is (isometrically isomorphic to) a closed
subalgebra of  ${\mathcal A}^{**}$, the pair ($\mathcal A$,
${\mathcal A}^{**}$) is  an admissible pair, and this
provides a short proof for Proposition~\ref{M1} (by virtue of
Theorem~\ref{M2}).

(ii)  Lau and Loy in their extensive work \cite{LL} gave a train of
admissible pairs, in particular in the group algebras
$L^{\infty}(G)$, $VN(G)$ and $PM_p(G)$. For instance, we quote some
of them which are related to $L^{\infty}(G)$, whose details can be
found in \cite[Sections 3 and 4]{LL}.

$\bullet$ $(M(G/N), M(G))$, where $N$ is a compact normal subgroup of $G$.

$\bullet$  $(M(G), {\mathcal X}^*)$, for every introverted subspace
${\mathcal X}$ of $L^\infty(G)$ with $C_0(G)\subseteq {\mathcal
X}\subseteq C(G)$. In particular, the pairs $(M(G), LUC(G)^*)$ and  $(M(G), WAP(G)^*)$.

$\bullet$ $(AP(G)^*, WAP(G)^*)$.

$\bullet$ $(WAP(G/N)^*, WAP(G)^*)$, where $N$ is a closed normal subgroup of $G$.

$\bullet$ The pairs  $({\mathcal X}^*, L^1(G)^{**})$ and
$(L^1(G/N)^{**}, {L^1(G)}^{**})$, in which $G$ is amenable, $N$ is a
closed normal subgroup of $G$ and $\mathcal X$ is a non-zero
$w^*-$closed self-adjoint translation invariant subalgebra of
$L^{\infty}(G)$.

One may also find some other admissible pairs related to  $VN(G)$ and $PM_p(G)$ in Sections 6 and 7 of \cite{LL}.

According to Theorem~\ref{M2}, in each of the above mentioned pairs,
the $n$-weak amenability of the second component implies that of the
first.
\end{remark}

\end{document}